\documentclass{llncs}
\usepackage{makeidx}  % allows for indexgeneration

\usepackage{epsfig,latexsym,amsfonts,amsmath,amssymb,color,mathrsfs,dsfont}
\usepackage{graphicx}
\usepackage{tikz}
\usepackage{mathtools}

\def\dvg{{\rm div}}

\def\IntO{\int\limits_{\Omega}}

\newcommand\be{\begin{eqnarray*}}
\newcommand\ee{\end{eqnarray*}}
\newcommand\ben{\begin{eqnarray}}
\newcommand\een{\end{eqnarray}}

\def\W{{\sf W}}

\def\benn{\begin{eqnarray*}}
\def\eenn{\end{eqnarray*}}
\def\BR{{\mathbb R}}

\newcommand\dx{\, dx}

\newcommand{\om}{-}
\newcommand{\op}{+}

\newcommand{\Oum}{\Omega ^u_\om}
\newcommand{\Oup}{\Omega ^u_\op}
\newcommand{\Oun}{\Omega ^u_0}

\newcommand{\Ovm}{\Omega ^v_\om}
\newcommand{\Ovp}{\Omega ^v_\op}
\newcommand{\Ovn}{\Omega ^v_0}

\newcommand{\operm}[1]{(#1)_\om}
\newcommand{\operp}[1]{(#1)_\op}

\begin{document}
\frontmatter          % for the preliminaries
\pagestyle{headings}  % switches on printing of running heads
\addtocmark{Verification of identities} % additional mark in the TOC

\mainmatter              % start of the contributions
\title{Verifications of primal energy identities for variational problems with obstacles}
\titlerunning{Verification of energy identities}  % abbreviated title (for running head)
\author{Sergey Repin \inst{1} \and Jan Valdman\inst{2}}
\authorrunning{Sergey Repin et al.} % abbreviated author list (for running head)
\tocauthor{Sergey Repin and Jan Valdman}
\institute{V.A.Steklov Institute of Mathematics in St.-Petersburg, 191011, Fontanka 27, Sankt--Petersburg, Russia  and University of Jyv\"askyl\"a, P.O.Box 35, FI-40014,  Finland
\and
Institute of Mathematics and Biomathematics, Faculty of Science, University of South Bohemia, Brani\v sovsk\' a 31, \v{C}esk\'e Bud\v ejovice, , Czech Republic, CZ--37005 and Institute of Information Theory and Automation, Academy of Sciences, Pod vod\'{a}renskou v\v{e}\v{z}\'{\i}~4, CZ--18208~Praha~8, Czech Republic
\email{jvaldman@prf.jcu.cz}}

\maketitle              % typeset the title of the contribution

\begin{abstract}
We discuss error identities for two classes of free boundary problems generated by obstacles. The identities suggest true forms of the respective error measures which consist of two parts: standard energy norm and a certain nonlinear measure. The latter measure controls (in a weak sense)
approximation of free boundaries. Numerical tests confirm sharpness of error identities and show that in different examples one or another part of the error measure may be dominant.
\keywords{variational problems with obstacles, coincidence set, error identities}
\end{abstract}

\section{Introduction}
%Variational inequalities form an important class of nonlinear
%models that describe free boundary phenomena arising in various
%applied problems (see, e.g., G.~Duvaut and J.~L. Lions \cite{DuLi} and other publications cited therein). Usually free boundaries separate regions where  solutions possess quite different physical properties. Therefore, any reliable information on the shape and location of such a boundary is very important.
%Qualitative properties of free boundaries are studied by purely analytical (a priori) methods unlike quantitative
%information, which in the vast majority of cases can be obtained only by computational methods.
%In this context, it is necessary to know which quantitative information could be indeed extracted
%from a numerical solution.
New types of error identities were recently derived \cite{ReVa2017} for two types of  inequalities generated by obstacle type conditions: a classical obstacle problem and a two-phase obstacle problem. Both problems belong to the class of variational problems
\ben
\label{1.1}
\inf\limits_{v\in V}J(v),\qquad J(v)=G(\Lambda v)+F(v),
\een
where $\Lambda:V \rightarrow Y$ is a bounded linear operator,
 $ G:Y\rightarrow \BR$ is a convex, coercive, and lower semicontinuous functional,
 $F:V\rightarrow\BR$ is another convex lower semicontinuous functional, and
 $Y$ and $V$ are reflexive Banach spaces. Henceforth, we use results of \cite{Repin2000} related to 
derivation of a posteriori error estimates for this class of problems.

\subsection{The classical obstacle problem}

The classical obstacle problem (see, e.g. \cite{DuLi,EkTe}) is characterized by 
\be
%\Lambda v=\nabla v,\qquad 
G(\Lambda v)=\frac12 \IntO  A \nabla v\cdot\nabla v \dx, \qquad F(v)=-\IntO fv\,dx +\chi_K(v),
\ee
where the characteristic functional  is defined as 
\be
\chi_K(v):=\left\{
\begin{array}{cc}
0\;&{\rm if}\;\phi\leq v\leq \psi,\\
+\infty& {\rm else}
\end{array}
\right.
\ee
and the admissible set reads 
\be
&&K:=\{v\in  V_0:=H^1_0(\Omega) \,\mid\,\phi(x)\,\leq v(x)\, \leq \psi(x)\;
{\rm a.e.\,in\,}\Omega \}.
\ee
Here, $H^1_0(\Omega) $ denotes the Sobolev space of functions vanishing
on $\partial\Omega$ (hence we consider the case $u_D=0$),
 $\Omega\subset\mathbb{R}^d$ ($d \in \{1, 2, 3\}$) is a bounded domain with a Lipschitz continuous boundary $\partial\Omega$ and 
$\phi,\psi\in H^2(\Omega )$ are two given functions (lower and upper obstacles) such that
\be
&&\phi(x)\leq 0\;{\rm on}\;\partial\Omega,\quad
\psi(x)\geq 0\;{\rm on}\;\partial\Omega,
\quad\phi(x)\leq\psi(x),\quad\forall x\in \Omega.
\ee
It is assumed that $A $ is a 
symmetric matrix subject to the condition
\ben \label{A_assumptions}
A(x) \xi\cdot\xi\geq c_1 \, |\xi|^2\qquad c_1>0,\qquad \forall \xi \in
\mathbb{R}^{d}
\een
almost everywhere in $\Omega $. Under the assumptions made, the
unique solution $u \in K$ exists.  The mechanical motivation of the obstacle problem is to find the equilibrium position of an elastic membrane whose boundary is held fixed, and which is constrained to lie between given lower and upper obstacles $\phi$ and $\psi$.

\subsection{The two-phase obstacle problem}
The functional $J(v)$  of the two-phase-obstacle problem (see, e.g. \cite{ShahgolianUraltsevaWeiss}) is defined by the relation 
\begin{equation}
J(v) := \IntO \Big( \frac{1}{2} A \nabla v\cdot\nabla v  - f v + \alpha_{\op} \operp{ v } + \alpha_{\om} \operm{ v }  \Big) \dx.    \label{eq:energy}
\end{equation}
The functional $J(v)$  is minimized on the set
$$
V_0+u_D:=\{v=v_0+u_D\,:\, v_0 \in V_0,\;
u_D \in H^1(\Omega)\}. 
$$
Here $u_D$ is a given bounded function that defines the boundary condition
($u_D$   may attain both positive and negative values on different parts of the boundary $\partial \Omega$). It is assumed that  the coefficients $\alpha_{\op},
\alpha_{\om}: \Omega \rightarrow \mathbb{R}$ are  positive
constants (without essential difficulties the consideration and main results can be extended to the case where they are positive
Lipschitz continuous functions). Also, it is assumed that
 $f \in L^{\infty}(\Omega)$, $A \in L^{\infty}(\Omega, \mathbb{R}^{d \times d})$, 
 and  the condition \eqref{A_assumptions} holds. 
 Since the functional $J(v)$ is strictly
convex and continuous  on $V$, existence and uniqueness of a minimizer $u \in V_0+u_D$ is guaranteed by well known results
of the calculus of variations. The mechanical motivation of the two-phase obstacle problem is to  find the equilibrium position of an elastic membrane in the two-phase matter with different gravitation densities related to $\alpha_{\om}$ and $\alpha_{\op}$.

\section{Error identities}
The solution $u$ of the classical obstacle problem divides $\Omega $ into three sets:
\ben
&\Oum:=&\{x\in\Omega \,\mid\, u(x)=\phi(x) \}\,, \nonumber  \\
&\Oup:=&\{x\in\Omega \,\mid\, u(x)=\psi(x) \}\,, \label{coincidence_sets}    \\
&\Oun:=&\{x\in\Omega\,\mid\,\phi(x)<u(x)<\psi(x)\}\,. \nonumber 
\een
The sets  $\Oum$ and
$\Oup$  are the {\em lower} and
{\em upper coincidence sets}
and  $\Oun$
is an open set, where $u$ satisfies the Poisson equation 
$\dvg  (A \nabla u) + f=0$. Thus, the problem involves {\em free boundaries}, which are unknown a priori. Let $v$ be an approximation of $u$. It defines approximate  sets
\ben 
&\Ovm:=&\{x\in\Omega \,\mid\, v(x)=\phi(x)\}\,,  \nonumber \\
&\Ovp:=&\{x\in\Omega \,\mid\, v(x)=\psi(x)\}\,, \label{coincidence_sets_approximate}\\
&\Ovn:=&\{x\in\Omega\,\mid\,\phi(x)<v(x)<\psi(x)\} \nonumber \,.
\een
Notice that unlike the sets in (\ref{coincidence_sets}), the sets 
(\ref{coincidence_sets_approximate}) are known.

\begin{theorem}[\cite{ReVa2017}] \label{theorem1}
Let $v \in K$ be any approximation of the exact solution $u \in K$ of the classical obstacle problem. Then it holds
\ben
\label{id:primal}
&\frac12 || \nabla (u - v) ||_A^2 + \mu_{\phi\psi}(v) = J(v) - J(u), \label{MM_final_obstacle} 
\label{id:dual}
\een 
where  
\ben
\mu_{\phi\psi}(v) :=\!\int\limits_{\Oum}\!\! \W_\phi (v-\phi)\dx+  \int\limits_{\Oup}\!\! \W_\psi(\psi-v)\dx,% =\!\!\!\! 
%\int\limits_{\Oum \setminus \Ovm }\!\!\!\! \W_\phi (v-\phi)\dx+  \!\!\!\!\int\limits_{\Oup \setminus \Ovp }\!\!\!\! \W_\psi(\psi-v)\dx,
\label{measure_primal_classical}
\een
and
$\W_\phi:=-(\dvg A\nabla \phi+f), \W_\psi:=\dvg A\nabla \psi+f$
are two nonnegative %(due to  \eqref{equilibrium_of_p*}) 
weight functions generated by the source term $f$, the obstacles $\psi, \phi$ and the diffusion $A$. 
\end{theorem}
Here, $\mu_{\phi\psi}(v)$ represents a certain (non-negative) measure, which
controls (in a weak integral sense) whether or not the function $v$
coincides with obstacles $\psi, \phi$ on true coincidence sets $\Oum$
and $\Oup$. 

\begin{remark}
The error identity \eqref{id:dual}  was derived for the homogeneous boundary condition $u=0$ on $\partial \Omega$, but it is possible to extend it in the same form to  for the nonhomogeneous boundary condition $u \not =0$ on $\partial \Omega$. 
\end{remark}

For the two-phase obstacle problem, we introduce two decompositions  (diffent from the classical obstacle problem)
of $\Omega$ associated with the minimizer $u$
and an approximation $v$:
\ben
&\Oum:=\{x\in \Omega \,\mid\, u(x)<0\}, \nonumber  \\
&\Oup:=\{x\in \Omega \,\mid\, u(x)>0\},  \label{coincidence_sets_two-phase}\\
&\Oun:=\{x\in \Omega \,\mid\, u(x)=0\}, \nonumber
\label{Omega_u_two-phase_obstacle}
\een
and
\ben 
&\Ovm:=\{x\in \Omega \,\mid\, v(x)<0\}, \nonumber\\
&\Ovp:=\{x\in \Omega \,\mid\, v(x)>0\}, \label{coincidence_sets_approximate_two-phase}\ \\
&\Ovn:=\{x\in \Omega \,\mid\, v(x)=0\}. \nonumber
\label{Omega_v_two-phase_obstacle}
\een
These decompositions generate exact and approximate free boundaries. If we introduce new sets
\be
 \omega_+:=\Ovp \cap \Oun, \quad
 \omega_-:=\Ovm \cap \Oun, \quad
\omega_\pm:=\left\{\Ovp \cap \Oum \right\}\cup \left\{\Ovm \cap \Oup \right\},
\ee
we can formulate an error identity for the two-phase obstacle problem.

\begin{theorem}[\cite{ReVa2015},\cite{ReVa2017}] \label{theorem2}
Let $v \in V_0+u_D$ be any approximation of the exact solution $u \in V_0+u_D$ of the two-phase obstacle problem. Then it holds
\ben
&&\frac12 || \nabla (u - v) ||_A^2 + \mu_\omega(v) = J(v) - J(u), \label{MM_final_obstacle_two-phase} 
\een
where
\ben 
\label{measure_primal_two-phase}
 \mu_\omega(v) :=\int\limits_{ \omega} \alpha(x)|v| \dx, \quad \omega:= \omega_+ \cup \omega_- \cup \omega_\pm
   %=\alpha_\op\int\limits_{\omega^0_+} v \dx
%- \alpha_\om\int\limits_{\omega^0_-} v \dx
%+ (\alpha_{\op} +\alpha_\om)\int\limits_{\omega_\pm }|v| \dx
\een
and
\ben
\alpha(x):=\left\{
\begin{array}{ll}
\alpha(x)=\alpha_\op\quad & {\rm if}\;x\in \,\omega_+,\\
\alpha(x)=\alpha_\om\quad & {\rm if}\;x\in \,\omega_-,\\
\alpha(x)=\alpha_\op + \alpha_\om \quad & {\rm if}\;x\in \,\omega_\pm.
\end{array}
\right.
\een
\end{theorem}
Here, $\mu_{\omega}(v)$ represents another nonlinear measure (which differs from $\mu_{\phi\psi}$).
%This measure is nonnegative and vanishes if and only if $\Ovp$ coincides with $\Oup$ and $\Ovm$ coincides with $\Oum$.

\section{Numerical verifications}
We verify a posteriori error identities \eqref{MM_final_obstacle} and \eqref{MM_final_obstacle_two-phase} for both obstacle problems and focus on interpretation of their nonlinear measures $\mu_{\phi\psi}(\cdot)$ and $\mu_{\omega}(\cdot)$. Another goal is to present examples with
different balance between two components of the overall error measure.

\subsection{The classical obstacle problem in 2D}
We assume a 2D example taken from \cite{NSV}. In this example,  
$\Omega=(-1,1)^2, A=\mathbb I, \phi=0, \psi=+\infty$
It is known that for 
\begin{equation*}
f(x,y)=\left\{
\begin{array}{lrl} -16(x^{2}+y^{2})+8R^2 & \quad\textrm{if} & \sqrt{x^2+y^2}>R\\
 -8(R^{4}+R^{2})+8R^{2}(x^{2}+y^{2}) & \quad\textrm{if} & \sqrt{x^2+y^2}\leq R
\end{array} \right. , 
\end{equation*}
where $R\in[0,1)$ is given, the exact solution to the obstacle problem reads
\begin{equation*}
u(x,y)=\left\{
\begin{array}{lrll}
\left(\max\{x^{2}+y^{2}-R^{2},0\}\right)^2 & \quad \mbox{if} &(x,y) \in \Omega \\
\left(x^{2}+y^{2}-R^{2}\right)^2 & \quad \mbox{if} &(x,y) \in \partial\Omega
\end{array} \right. .
\end{equation*}
The corresponding energy can be computed (see \cite{HaVa2})  and it reads
$$
J(u) = 192\left(\frac{12}{35}-\frac{28R^{2}}{45}+\frac{R^{4}}{3}\right)
-32R^{2}\left(\frac{28}{45}-\frac{4R^{2}}{3}+R^{4}\right)
+\frac{2}{3}\pi R^{8}.
$$ 

\begin{figure}
  \centering
%\begin{minipage}{5.2cm}
%\includegraphics[width=\textwidth]{perturbation_n_0} \\
%\end{minipage}
%\begin{minipage}{1cm}
%\end{minipage}
%\hspace{1cm}
%\begin{minipage}{5.2cm}
 %\includegraphics[width=\textwidth]{coincidence_set_n_0} 
%\end{minipage}
%$$\mbox{Perturbation function for } r=0.4, k=0 \mbox{ and the corresponding coincidence set } \Ovm. $$
%
%\begin{minipage}{5.2cm}
%\includegraphics[width=\textwidth]{perturbation_n_8} \\
%\end{minipage}
%\begin{minipage}{1cm}
%\end{minipage}
%\hspace{1cm}
%\begin{minipage}{5.2cm}
 %\includegraphics[width=\textwidth]{coincidence_set_n_8} 
%\end{minipage}
%$$\mbox{Perturbation function for } r=0.3, k=8 \mbox{ and the corresponding coincidence set } \Ovm. $$
\begin{minipage}{5.2cm}
\includegraphics[width=\textwidth]{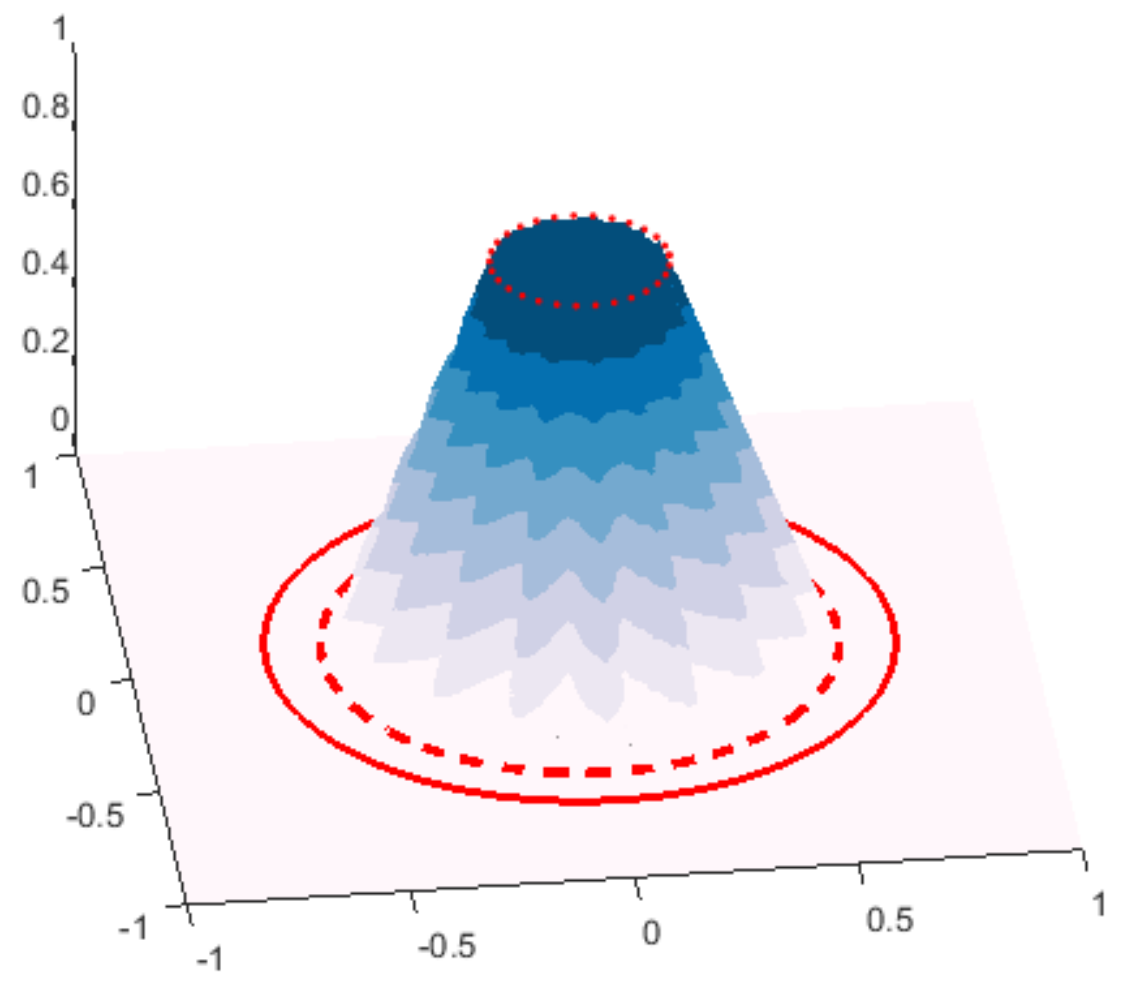} 
\end{minipage}
\begin{minipage}{1cm}
\end{minipage}
\hspace{1cm}
\begin{minipage}{5.2cm}
 \includegraphics[width=\textwidth]{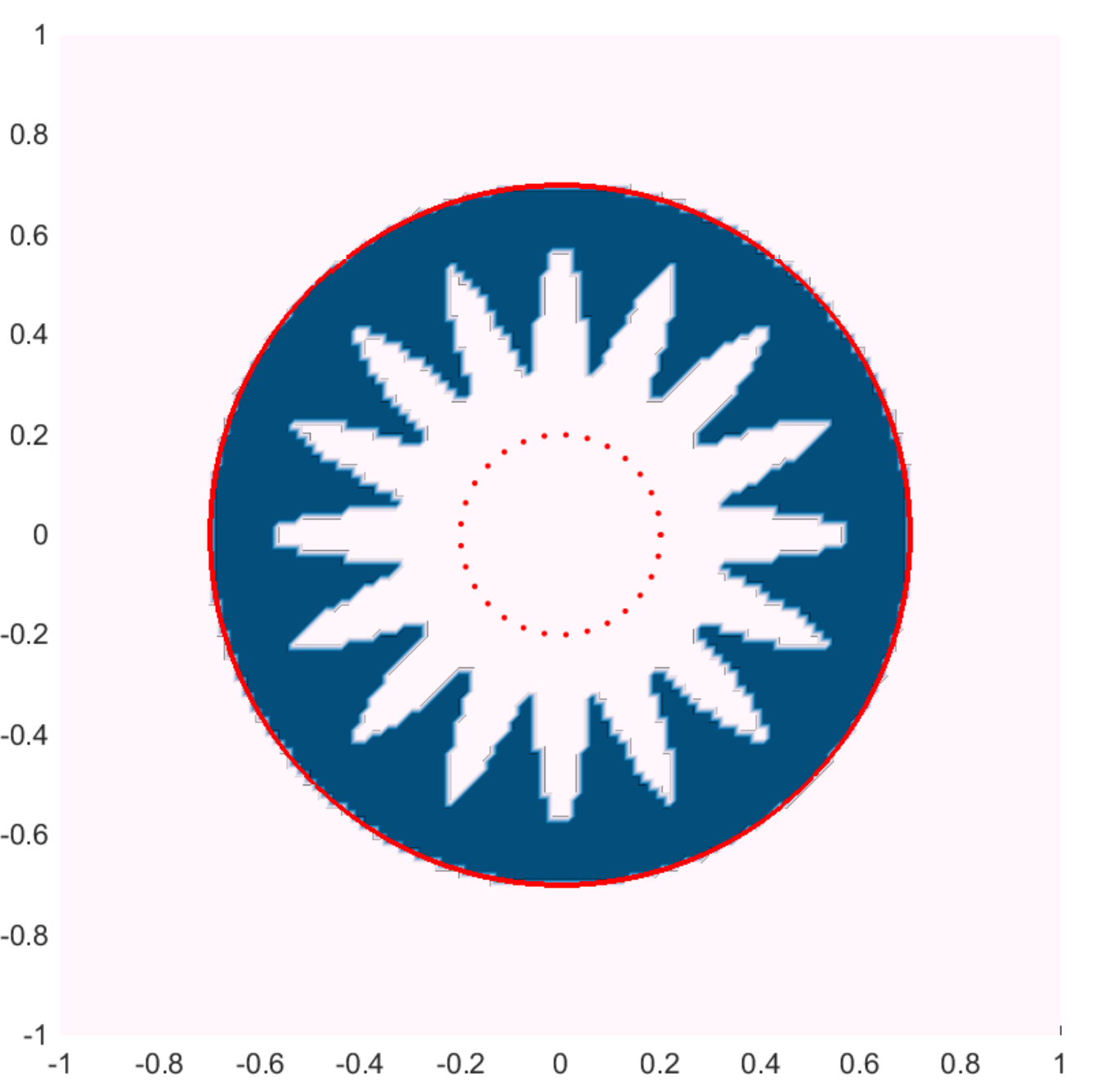} 
 %coincidence_set_n_16
\end{minipage}
\caption{A perturbation function $w$ (left) generated 
by parameters $R=0.7, r=0.2, k=16$ and the same corresponding coincidence set $\Ovm$ (right) for all approximative solutions $v = u + \epsilon \, w$, where 
$\epsilon >0$.  The boundary of $\Oum$ is indicated by the full circle, the inner radius  $r$ by the dotted circle and the intermediate radius $\frac{r+3 R}{4}$ by the dashed circle. }
\label{fig2}  
\end{figure}

We consider approximations $v$ in the form
\begin{equation}
v_\epsilon: = u + \epsilon \, w, \label{v_2D} 
\end{equation}
where $\epsilon > 0$ is a given amplitude and $w$ is a solution perturbation defined in polar coordinates $(\rho, \theta)$ as 
\begin{equation}
w(\rho,\theta): = \left\{\begin{array}{ll}
        1, & \text{if } \rho \leq r \\
        1 - \frac{\rho - r}{\tilde{r}(\theta) - r}, & \text{if } r \leq \rho \leq \tilde{r}(\theta) \\
        0, & \text{if } \rho \geq \tilde{r}(\theta) 
        \end{array}\right. .
\end{equation}
Here, $0 < r < R$ is given internal radius and a variable radius $\tilde{r}(\theta)$ is defined as 
\begin{equation}
 \tilde{r}(\theta):= r + (R - r) \left( \frac{ 2 + \cos (k \theta )}{4} \right) \label{var_radius}
\end{equation}
for some $k \in \mathbb{Z}$.
 This construction ensures that 
\begin{equation}
r < \frac{3 r+R}{4} \leq \tilde{r}(\theta) \leq \frac{r+3 R}{4} < R \label{bounds_radius}
\end{equation}
and consequently $\nabla w$ is bounded. An examples of perturbations $w$ is visualized in Figure \ref{fig2} together with corresponding coincidence sets $\Ovm$. For given $k$ and $r$, there is always a convergence in the energy error
\ben v_\epsilon \rightarrow u \qquad (\mbox{in } K )  \qquad  \mbox{ as } \epsilon \rightarrow 0 \label{convergence_energy}\een
and consequently the nonlinear measure must also converge
\ben \mu_{\phi\psi}(v_\epsilon)   \rightarrow  \mu_{\phi\psi}(u)  =0 \qquad  \mbox{ as } \epsilon \rightarrow 0. \label{convergence_nonlinear_measure}\een
 It should be noted the shape of $\Omega_{-}^{v_\epsilon}$  depends on $k$ and $r$ only and it is completely independent of $\epsilon$. Therefore, 
%the approximate coincidence set 
$\Omega_{-}^{v_\epsilon}$ never approximates 
%the exact coincidence set 
$\Oum = \{ x \in \Omega: ||x|| \leq R \}$  for any choice of $\epsilon$! 
\begin{table}
\begin{center}
\begin{small}
\begin{tabular}{|l|c|c|c|c|r|}
\hline
$\epsilon$  & $\frac{1}{2}\|\nabla(u-v_{\epsilon})\|^{2}_{A}$ & $\mu_{\phi\psi}(v_{\epsilon})$ & $J(v_{\epsilon})-J(u)$  & $\kappa(v_{\epsilon}) \, [\%]$ \\
\hline
1.0000  & 7.1531e+00 & 4.4311e+00 & 1.1584e+01 & 38.2512 \\
0.1000  & 7.1531e-02 & 4.4311e-01 & 5.1464e-01  & 86.1008 \\
0.0100  & 7.1531e-04 & 4.4311e-02 & 4.5027e-02 & 98.4113 \\
0.0010  & 7.1531e-06 & 4.4311e-03 & 4.4388e-03 & 99.8388 \\
0.0001  & 7.1531e-08 & 4.4311e-04 & 4.4375e-04  & 99.9839 \\
\hline 
\end{tabular}
\caption{The error identity parts computed for various $v_{\epsilon}=u+\epsilon w$, where the exact  coincidence set $\Oum$ is represented by the circle of the radius $R=0.7$ and the perturbation $w$ is defined by the choice $r=0.2, k=16$.} 
\label{table2}
\end{small}
\end{center}
\vspace{-1cm}
\end{table}

Table \ref{table2} reports on values of terms
%$$\frac{1}{2}\|\nabla(u-v_{\epsilon})\|^{2}_{\Omega,A}, \quad \mu_{\phi\psi}(v_{\epsilon}), \quad J(v_{\epsilon})-J(u)$$ 
in the energy identity \eqref{MM_final_obstacle} for few approximations $v_{\epsilon}$, where $\epsilon$ decreases to $0$ and $u$ and $w$ are given by the choice of $R$ and $r,k$. If $\epsilon$ tends to zero, the term $\frac{1}{2}\|\nabla(u-v_{\epsilon})\|^{2}_{A}$ converges quadratically to 0 and the nonlinear measure $\mu_{\phi\psi}(v_{\epsilon})$ only linearly to 0. The contribution of the nonlinear measure to the energy identity is measured by the quantity
\ben \label{ratio}
\kappa(v_{\epsilon}) := 100 \, \frac{\mu_{\phi\psi}(v_{\epsilon})}{J(v_{\epsilon})-J(u)} \quad [\%].
\een
We see in this example, the contribution of $\mu_{\phi\psi}(v_{\epsilon})$ dominates over the contribution of $\frac{1}{2}\|\nabla(u-v_{\epsilon})\|^{2}_{A}$. 

\subsection{The two-phase obstacle problem in 1D}
This subsection extends results of \cite{ReVa2015}. We consider the two-phase obstacle problem in 1D from \cite{Bo}. Here, $\Omega=(-1, 1), f=0, A=\mathbb I, \alpha_{\oplus}=\alpha_{\ominus}=8$
and the Dirichlet boundary conditions $u(-1)=-1, u(1)=1.$ The exact solution is given by 
$$ u(x)=\left\{\begin{array}{ll} -4x^2-4x-1 , \quad & x \in \left[-1, -0.5 \right], \\
                                 0, \quad & x \in \left[-0.5, 0.5 \right], \\
                                 4x^2-4x+1, \quad & x \in \left[ 0.5, 1 \right]
                                    \end{array} \right.
                                                           $$
and $J(u)=5\frac{1}{3}.$ 
%The function $u$ splits the domain $\Omega$ in subdomains
%$$ \Omega^-_u=(-1,-0.5), \qquad \Omega^0_u=\left< -0.5, 0.5 \right>, \qquad \Omega^+_u=(0.5, 1).$$
%The error identity \eqref{mainidentity} is valid for any approximation $v \in H^1(-1,1)$ satisfying boundary conditions $v(-1)=-1, v(1)=1$. 
\begin{figure}[t]
  \centering
\begin{minipage}{12cm}
	%$$\mbox{ The approximation $v=v_N$ is the interpolant of $u$ with $N=2$ nodes:}$$
  %\includegraphics[width=\textwidth]{two-phase_obstacle_1D_N_2} \\
	%%\vspace{-0.2cm}
%
		%$$\mbox{ The approximation $v=v_N$ is the interpolant of $u$ with $N=5$ nodes:}$$
  %\includegraphics[width=\textwidth]{two-phase_obstacle_1D_N_5} \\
%			$$\mbox{ The approximation $v=v_N$ is the interpolant of $u$ with $N=6$ nodes:}$$
  \includegraphics[width=\linewidth]{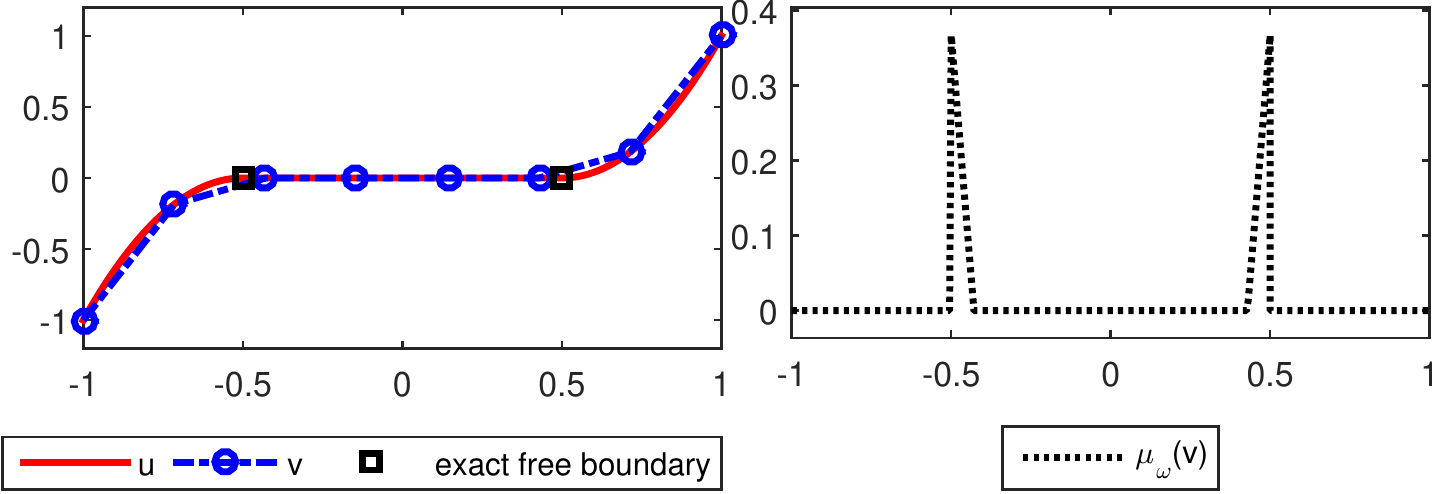}
  % height=0.4\linewidth
\end{minipage}
  \caption{Exact solution $u$ of the two-phase obstacle problem and its approximations $v_N$ (left) for $N=8$ and the distributions of $\mu_{\omega}(v_N)$ (right).}
\label{fig1}
\end{figure}
We consider a sequence of approximations
 $$v_N(x)=I_N(u)(x), \quad x \in \left[ -1,1 \right],$$
where $I_N$ (for $N=2, 3, \dots $) denotes a piecewise linear nodal interpolant of the function $u$ in $N$ uniformly distributed nodes 
$\{ -1, -1+h, \dots, 1-h, 1\},$ 
where $h=2/(N-1)$. 
%The simplest choice $N=2$ (and $N=3$ as well) generates an approximation $v_{2}(x)=x, x \in \left< -1,1 \right>$ and the direct calculation \cite{ReVa2015} shows 
%that $J(v_{2})=9$ 
%$$\frac12\|\nabla(u-v_{2})\|^{2}_{\Omega,A}=1\frac{2}{3}, \qquad \mu_{\omega} (v_{2})=2,  \qquad J(v_2)-J(u)=3\frac{2}{3},$$
%so that the energy identity \eqref{MM_final_obstacle_two-phase} is satisfied. 
%Let us also consider an approximation $v_5$. %depicted in Figure \ref{fig1}. 
%Two of its five interpolation nodes lie on the exact free boundary at $x=\pm 0.5$ and sets $\Omega^-_{v_5}, \Omega^0_{v_5}, \Omega^+_{v_5}$  
%$$ \Omega^-_{v_5}=(-1,-0.5), \qquad \Omega^0_{v_5}=\left< -0.5, 0.5 \right>, \qquad \Omega^+_{v_5}=(0.5, 1)$$
%coincide with $\Omega^-_u, \Omega^0_u, \Omega^+_u$ respectively. The approximation $v_5$ clearly approximates the free boundary exactly but provides only very coarse approximation of $u$ on $\Omega^-_v$ and $\Omega^+_v$. 
Table \ref{table1} reports on terms in the energy identity \eqref{MM_final_obstacle_two-phase}
%$$\frac{1}{2}\|\nabla(u-v_N)\|^{2}_{\Omega,A}, \quad \mu_{\omega}(v_N), \quad J(v_N)-J(u)$$ 
for some increasing values of $N$. In general, it holds
$$\mu_{\omega}(v_N)=0 \qquad \mbox{for } N=4k+1, k \in \mathbb{N}. $$
In these cases, two interpolation nodes lie on the exact free boundary at $x=\pm 0.5$ 
and sets $\Omega^-_{v_N}, \Omega^0_{v_N}, \Omega^+_{v_N}$  
coincide with $\Omega^-_u, \Omega^0_u, \Omega^+_u$. For all other approximations $v_N$ (see Figure \ref{fig1} for $N=8$), it holds $\mu_{\omega}(v_N) > 0$. The contribution of the nonlinear measure to the energy identity is measured by the quantity
\ben \label{ratio_two-phase}
\kappa(v_N) := 100 \, \frac{\mu_{\omega}(v_N)}{J(v_N)-J(u)} \quad [\%].
\een
We see in this benchmark, the contribution of $\frac{1}{2}\|\nabla(u-v_N)\|^{2}_{A}$ dominates over the contribution of the nonlinear measure term $\mu_{\omega}(v_N)$.

%\begin{remark}
%Terms of Table \ref{table1} are evaluated by the trapezoidal quadrature on 1 million uniformly %uniformly distributed subintervals.
% A MATLAB software 
%extends the original code of \cite{ReVa2015} and 
%it is available for own testing \cite{code_two-phaseobstacle_1D}. 
%\end{remark}

\begin{table}[h]
\begin{center}
\begin{small}
\begin{tabular}{|l|c|c|c|c|c|}
\hline
$N$   & $\frac{1}{2}\|\nabla(u-v_N)\|^{2}_{A}$ & $\mu_{\omega}(v_N)$ & $J(v_N)-J(u)$  & $\kappa(v_{N}) \, [\%]$ \\
\hline
2  & 1.67e+00 & 2.00e+00 & 3.67e+00 & 54.55 \\ 
5  & 6.67e-01 & 00 & 6.67e-01 & 0.00 \\ 
6  & 3.59e-01 & 7.20e-02 & 4.31e-01 & 16.72 \\ 
7  & 2.59e-01 & 7.41e-02 & 3.33e-01 & 22.22 \\ 
8  & 2.16e-01 & 2.62e-02 & 2.42e-01 & 10.82 \\ 
9  & 1.67e-01 & 00 & 1.67e-01 & 0.00 \\ 
10  & 1.20e-01 & 1.23e-02 & 1.32e-01 & 9.33 \\ 
%11  & 9.87e-02 & 1.60e-02 & 1.15e-01 & 13.95 \\ 
%12  & 8.78e-02 & 6.76e-03 & 9.45e-02 & 7.15 \\ 
%13  & 7.41e-02 & 00 & 7.41e-02 & 0.00 \\ 
%14  & 5.92e-02 & 4.10e-03 & 6.33e-02 & 6.47 \\ 
%15  & 5.15e-02 & 5.83e-03 & 5.73e-02 & 10.17 \\ 
30  & 1.23e-02 & 3.69e-04 & 1.27e-02 & 2.91 \\ 
60  & 3.06e-03 & 4.38e-05 & 3.10e-03 & 1.41 \\ 
120  & 7.53e-04 & 5.34e-06 & 7.58e-04 & 0.70 \\ 
\hline
\end{tabular}
\caption{The error identity terms computed for various approximation $v_N$.}
\label{table1}
\end{small}
\end{center}
\vspace{-1cm}
\end{table}

\subsubsection*{Acknowledgments.} 
The first author acknowledges the support of RICAM during Special Semester on Computational Methods in Science and Engineering 2016, Linz, Austria. The
second author has been supported by GA CR through the projects GF16-34894L and 17-04301S.

\end{document}